\documentstyle{article}
\newtheorem{theo}{\bf Theorem}[section]

\newtheorem{defi}{\bf Definition}[section]
\newtheorem{nota}{\bf Notation}[section]
\newtheorem{lem}{\bf Lemma}[section]

\begin{document}

\begin{center}
{\bf {\Large Orthogonal Main Effect Plans on blocks of small size.}}
\end{center}
\vskip5pt

\begin{center}
  Sunanda Bagchi,
   Theoretical Statistics and Mathematics Unit,
Indian Statistical Institute, Bangalore 560059, India.
\end{center}
\vskip10pt

{\bf Abstract}

In this paper we  define  the concept of orthogonality between two
factors "through another factor". Exploiting this property we have
been able to obtain orthogonal main effect plans (OMEP) on
non-orthogonal blocks requiring considerably smaller number of blocks than 
the existing methods.

We have also constructed  saturated partially orthogonal main effect plans 
(MEPs) for (i) an $n^4.2^3$ experiment and (ii) an $n^4.2.3$ experiment 
both on $4n$ runs. Here $n$ is an integer $\geq 3, n \neq 4 $.

 As particular cases, we have been able to accomodate four six-level factors 
on 8 blocks of size 4 each using the first method and on
24 runs using the second. 
 


\newpage
\section{Introduction}

Main effect plans (MEPs) with mixed levels are often required
for industrial experiments. The orthogonal main 
effect plans (OMEP) are, of course, the best option. However, due
to the divisibility conditions, an OMEP for an asymmetrical 
experiments often require a large number of runs. For this reason,
considerable attention has been devoted in the recent past to find 
main effect plans with the orthogonality condition relaxed to an
extent. Such plans were first proposed in Wang and Wu (1992). Subsequently many
others like Nguyen (1996), Ma, Fang and Liski (2000), Huang, Wu and Yen (2002) 
and Xu (2002) proposed and studied what they term "nearly orhogonal" arrays,
concentrating mostly on two- or three-level factors. In many of these plans, 
there are factors non-orthogonal to
three or more factors. As a result, in spite of the elegant combinatorial
properties, the precision of the estimates go down. 
 
   In our earlier paper [Bagchi (2006)] we have presented
"inter-class orthogonal" MEPs, where a factor is non-orthogonal
only to factors in its own class - to at most two factors in many of the 
plans. In this paper we continue the search for efficient MEPs
with run size not too big.    

Study of an orthogonal main effect plan with possibly non-orthogonal 
blocking was initiated in Mukerjee, Dey and Chatterjee (2001). They derived 
sufficient conditions for 
an OMEP to be universally optimal and also suggested a 
construction procedure for obtaining optimal OMEPs. In Bose and
Bagchi (2007) we came up with a new set of suficient conditions for OMEP 
on non-orthogonal blocks of size two, requiring smaller
number of blocks than Mukerjee, Dey and Chatterjee (2001). In this paper we 
generalise this idea. We define the concept of orthogonality between two
factors "through another factor" (which may or may not be a blocking factor).
Using this concept we have constructed an OMEP for a $4.n^3$ experiment on 
$n$ blocks of size 4 each, for $n \geq 5$. Treating the block factor as 
another treatment factor, one would obtain a $4.n^4$ experiment on $4n$ runs.
Although the main effects of this additional factor are not estimated with high
precision, [see (\ref {C-D})], there is considerable amount of reduction in 
number of runs. [We may recall that an OMEP for an $n^4$ experiment requires 
at least $n^2$ runs]. The four-level factor may also be replaced
by three two-level factors or one three and one two-level factor. 

We have presented another method showing that an OMEP for an $s^m$ experiment 
on $bk$ non-orthogonal blocks of size $ k \geq 2$ exists provided a connected 
binary block design with parameters $(b,k,s)$ exist. Here $m$ is the
maximum number of constraints of an orthogonal array of strength
two, $k^2$ runs and $k$ symbols. 
 In particular, we have obtained an OMEP for a $6^4$
 experiment on 8 blocks of size 4 each.

 \section{ An alternative to an existing MEP.}

We begin with a new main effect plan for a $ 3^4 2^3$ experiment
and compare it with the existing plan. Before that we need a few
notation.

\begin{nota} \label{MEP} A main effect (MEP) plan with a set $F$
of $m$ factors $ \{A,B \cdots \}$ and $n$ runs will be represented by an 
array $\rho(n,m; s_A \times s_B \cdots )$, termed {\bf mixed array}.
In $\rho$  rows will represent factors (in natural order) and
$S_A$ (respectively $s_A$) will denote the set (respectively number) of 
levels of the factor $A \in F$.  Often, different rows have the same number of 
symbols and we may denote the plan as $\rho(n,m; \prod_{i=1}^{k} 
(s_{i})^{m_i})$,  where $\sum_{i=1}^{k} m_{i} = m$.

The {\bf vector of unknown effects} of the levels of factor $A$ will be 
denoted by the $s_A \times 1$ vector {\bf $\alpha$}.
 The {\bf replication vector} of the factor A is the vector of
replication numbers of its levels in the natural order and will be
denoted by  {\bf$r^A$}. {\bf $R^A$} will denote the diagonal matrix with entries
as that of $r^A$ in the same order. 
{\bf $C_A$} will denote the {\bf coefficient matrix (C-matrix)} in the 
reduced normal equation obtained by eliminating the effects of all  
factors other than $A$.

 For two factors A,B, the {\bf incidence matrix} $N^{A,B}$ is an $s_A \times 
s_B$ matrix with the $(k,l)$th entry as the number of occurences
 of the level $k$ of $A$ and $l$ of $B$ together. $N^{A,B}= R^A$, when  $B=A$.
\end{nota}

For a plan represented by a mixed array $\rho(n,m; s_A \times s_B \cdots )$
let $Y$ denote the vector of yields. Then, our model is 

\begin{equation}
E(Y) = \mu 1_n + \sum \limits_{A \in F} X_A \alpha , \end{equation}

where $\mu$ is the general effect, $\alpha$ is as stated above
and $X_A$ is a 0-1 matrix as described below.
 For $A \in F$, the $(i,j)$the entry of $X_A$ is 1 if the ith
column of $\rho$ contains $j \in S_A$ in the row coresponding to
factor $A$ and 0 otherwise.

We present the well-known definition of an orthogonal array  for the sake of 
completeness.

\begin{defi}  An orthogonal array $OA(n,m,s_1\times \cdots \times s_m,t)$, 
having $m(\geq 2)$ rows, $n$ columns,  $s_1,\ldots ,s_m(\geq 2)$
symbols and strength $t(\leq m)$, is an $m\times n$ array, with elements in 
the $i$th row from a set of $s_i$ distinct symbols $(1\leq i\leq m)$, in 
which all $t$-tuples of symbols appear equally often (say
$\lambda$ times) as columns in every $t \times n$ subarray. \end{defi}
If $s_i=s, 1\leq i\leq m$, the OA is denoted by $OA(n,m,s,t)$. 
$\lambda$ is called the index of the OA. 
\vspace{.5em}
Let us recall the nearly orthogonal array $L_{12}' (3^42^3)$ of Wang and Wu 
(1992), referred to here as $A_{WW}(12)$. 

We shall suggest an MEP $A_1 (12)$  for a $3^4 2^3$ experiment
described as follows.

\begin{equation} \label{2^3array} \mbox{ Let } 
U_1 = \left [ \begin{array}{cccc}
                                     0 & 0 & 1 & 1 \\
                                     0 & 1 & 0 & 1 \\
                                     0 & 1 & 1 & 0 \\
\end{array} \right ].\\ \end{equation}

Let $V$ denote the array obtained by adding a row of all zeros
to $U_1$. Then,

$$ A_1 (12) = \left [ \begin{array}{ccc}
    V   &  V+1 & V + 2 \\
    U_1 & U_1 & U_1  \\ \end{array} \right ].\\$$

[Here addition is modulo 3].
 
Comparing the MEPs $A_{WW} (12)$ and $A_1 (12)$, we find the following.

1. In both the plans, each two-level factor is orthogonal to
every other factor.

2. In both the plans, no pair of three-level factors satisfy
proportional frequency condition.

3. The C-matrices of the three-level factors for the plan
$A_{WW} (12)$ are as follows.

\begin{equation} \label{C-ABC} C_Q = (7/3) K_3, \; Q=A,B,C,D. 
\end{equation}

 [Here $K_n = I_n - (1/n) J_n$, $J$ is the matrix of all-ones.]

Those for $A_N (12)$ are as follows.

$$C_Q = 3K_3, Q=A,B,C; \; Q_D = (4/3) K_3.$$

 We observe that $A_1 (12)$ provides more information to factors A,B,C 
than $A_{WW} (12)$ but less to D. Further, $\sum\limits_{Q=A,B,C,D}
C_Q$ is bigger for $A_1 (12)$, so that the total information
is more for this plan. Thus, $A_1 (12)$ is useful in situations
where all the factors are not equally important.

Now, we go to the deeper question. We note that in $A_{WW} (12)$,
for each pair $P \neq Q, P,Q \in \{A,B,C,D \}$, $N_{PQ}$ is the
incidence matrix of a balanced block design (BBD), the "best"
possible incidence matrix. On the other hand, for $A_1 (12)$,
neither of $N_{PQ}, P \neq Q$ is a "good" incidence matrix.
How is then, $C_A$, for instance, is bigger for $A_1 (12)$ ?

Before going to the investigation of this mystry,  we present
two plans which are obtained by  modifying $A_1 (12)$ slightly.

\begin{equation} \label{U2U3}
  \mbox{ Let } U_2 = \left [ \begin{array}{cccc}
  0 & 0 & 0 & 1 \\
  0 & 1 & 2 & 0  \\ \end{array} \right ] \\ 
\mbox{ and }
U_3 = \left [ \begin{array}{cccc} 0 & 1 & 2 & 3
\end{array} \right ].\\ \end{equation}

Then,
$$ A_i (12) =  \left [ \begin{array}{ccc}
    V   &  V+1 & V + 2 \\
    U_i & U_i  & U_i  \\ \end{array} \right ], i=2,3.\\$$

$ A_2 (12)$ and $ A_3 (12)$ are an MEPs for a $3^5.2$ and a $3^4.4$
experiment. The C-matrices of the first four (three-level) factors are, 
of course, same as
those for $ A_1 (12)$. The C-matrices of the new three-level factor and
the two-level factor of $ A_2 (12)$ are $3K_3$ and $3K_2$ respectively.
The C-matrix of the four-level factor in $A_3 (12)$ is $3K_4$.
[$K_n$ is as in the statement below (\ref {C-ABC}). 

To go into the mystry of "bigger" C-matrices of A,B,C in $A_i (12)$, we need 
some more  notation.

\begin{nota} \label{C-UV}For $L \subseteq F; U, V \in F\setminus L$, 
$C_{U,V;L}$ will denote the $(U,V)$th submatrix of the coefficient matrix 
(C-matrix) in the reduced normal equation  obtained by eliminating 
the effects of the member(s) of $L$. $C_{U,U;L}$ will be denoted by simply 
$C_{U;L} $. When $ L = S \setminus U$, $C_{U;L}$ will be denoted simply by 
$C_U$, to be consistent with Notation \ref {MEP}.
\end{nota}

\begin{defi} \label{orth-C} If in a mixed array the factors A,B, and C
satisfies the following condition, then we say that factors
A and B are mutually {\bf orthogonal through C}. 

\begin{equation} \label{orth-cond}
 N^{A,B} = N^{A,C} (R^C)^{-1} N^{C,B}\end{equation}
\end{defi}

The following lemma is immediate from the definition above.
\begin{lem} \label{orth-relative}
If in a mixed array a pair of factors A and B are orthogonal through C,
then $C_{A,B;C} =0$. \end{lem}

\begin{lem} \label{C-orth-relative} Suppose in a mixed array $\rho$ two 
factors A and  B satisfy the following condition.    

For every $Q \in \tilde{F} = F \setminus \{A,B\}$, A and Q are mutually  
orthogonal through B. Then, 

\begin{equation} C_A = C_{A;B}. \end{equation}
\end{lem} 

{\bf Proof :} We know that 

\begin{equation} \label{cmatcal} C_A = C_{A;B} - E_{A;B}(H_{A;B})^-(E_{A;B})^T,
\end{equation}

where $E_{A;B} = ((C_{A,Q;B}))_{Q \in \tilde{F}}$ and $H_{A;B} =
((C_{P,Q;B}))_{P \neq Q, \; P,Q \in \tilde{F}}$.

From the hypothesis, $E_{A:B} = 0$. Hence the result. $\Box$. 

Let us now look at the plans $A_i (12), i=1,2,3$. 
We observe the following property. 

\begin{lem} \label{C-mat for A(12)} For each of the plans $A_i (12), i=1,2,3$,  
 the C-matrices of factors A,B,C satisfy

\begin{equation} C_P = C_{P;D} , P \in \{A,B,C \}. \end{equation} \end{lem}

{\bf Proof :} 
It is easy to verify that each pair of factors $(P,Q), P,Q \in  \{A,B,C \}$ , are 
mutually orthogonal through $D$. Hence the
result follows from Lemma \ref {C-orth-relative}. $\Box$

Lemma \ref {C-mat for A(12)} is the clue of the mystry of the
"bigger" C-matrix in spite of the "poor" incidence matrices. We
note that for $A_{WW} (12)$  $E_{P:D} \neq 0$ for each $ P \in \{A,B,C \}$
and so, $C_P < C_{P;D}$. As a result, even though $C_{P;D}$ is bigger for 
$A_{WW} (12)$ than that for $A_1 (12)$,  $C_P$ becomes smaller.
   
  However, no such facility is available for the factor $D$. Further, the 
incidence matrices of D with the other factors are not "good". This is why $C_D$ is 
smaller (than that in $A_{WW} (12)$.   

Before going to the general construction in the next section, we present 
another example of a plan with two factors "orthogonal through another". 
The plan is for a $3^3$ experiment on 8 runs.

$A_8 = \left [ \begin{array}{cccccccc}
                                  0 & 0 &  0 & 0 & 1 & 1 & 2 & 2 \\
                                  0 & 2 &  0 & 2 & 0 & 1 & 0 & 1 \\
                                  2 & 0 &  0 & 2 & 1 & 0 & 0 & 1 \\
\end{array} \right ]$.

It is easy to verify that for $A_8$ the condition obtained
by interchanging $A$ and $C$ in (\ref {orth-cond}) is satisfied,
so that  factors B and C are orthogonal through A. As a result, 

\begin{equation} C_B = C_{B;A} =  
\left [ \begin{array}{ccc}
2  & -1 & -1 \\
-1 & 1  & 0 \\
-1 & 0  & 1 \\\end{array} \right ]. \end{equation}

We note that $C_B$ has spectrum $0^1.1^1.3^1$, which is only marginally
"smaller" than the spectrum $0^1.(9/4)^1.3^1$ of the C-matrix of
a three-level factor with the hypothetical OMEP with best
possible replication vector.

Factor $C$ also has the same C-matrix.

Since  $A$ is non-orthogonal to two others, $C_A$ is "smaller",
as shown below.

\begin{equation} C_A= (1/6)
\left [ \begin{array}{ccc}
4  & -2 & -2 \\
-2 & 7 & -5 \\
-2 & -5 & 7 \\\end{array} \right ].\end{equation}

The spectrum of $C_A$ is $0^1.1^12^1$.

\section{Two series of orthogonal main effect plans on blocks of
small size}

\begin{nota} \label {factbl}  An MEP for an $\prod_{i=1}^{m} s_{i}$ experiment laid out
on $b$ blocks of size $k_1.k_2, \cdots k_b$ each will be represented by a mixed
array $\rho(n,m+1; s_A \times s_B \cdots \times b$; the last row of $\rho$
represting the block factor. The set of factors will be
represented by $F \cup \{bl\}$, $bl$ denoting the block factor. 
$R^{bl}$ will denote the diagonal matrix with entries as $k_1.k_2, \cdots k_b$.
\end{nota}     

For the sake of easy reference, we write down the definition for
orthogonality through the block factor for an equi-block-sized plan
in a combinatorial form.

Let us consider an MEP on $b$ blocks of size $k$ each blocks represented by a 
mixed array $\rho(n,m+1; s_A \times s_B \cdots \times b)$. Fix two factors $A$ 
and $B$. For $ i \in S_A, j \in S_B$, let $U_{i,j}$ denote the number of 
columns in $\rho$ with A at level i and B at level j. Let 
$B_{i,j}$ denote the number of times the level combination (i,j) of A,B 
appear in the same block of $\rho$ (need not be in the same column). 

\begin{defi} \label{orth-bl} Consider a mixed array $\rho$ including a 
blocking factor of equal frequency $k$. Two factors A and B of $\rho$ are 
said to be mutually orthogonal through the block factor if 

\begin{equation} \label{orth-bleq} B_{i,j} = k U_{i,j}, \; \forall
i \in S_A, j \in S_B. \end{equation}
\end{defi}
 
 Note that $U_{i,j}$ is nothing but the $(i,j)$th entry of $N^{A,B}$, while 
$B_{i,j} = N^{A,bl} N^{bl,B}$. These together with the fact that $R^{bl} = 
k I_b$ implies that (\ref {orth-bleq}) is a special case of (\ref {orth-cond})
and hence Definition \ref {orth-bl} is a special case of 
Definition \ref {orth-C}.

{\bf Remark :} The case $k=2$ is considered in Theorem 3.1 of Bose and Bagchi 
(2007).

We now generalise the arrays $A_i(12), i=1,2,3$ to get an
infinite series of MEPs.

\begin{theo} \label{constr} (a) For an integer $n \geq 5$, there
exist saturated MEPs for the following asymmetrical experiments on $n$ 
blocks of size 4 each. 

(i) $n^3.2^3$, (ii) $n^3.2.3$ and (iii) $n^3.4$. 

(b) Plans (i) and (iii) are orthogonal while (ii) is almost orthogonal. 
Specifically, in all these plans, the two, three and four-level factors are 
orthogonal to each of the $n$-level factors as well as the block factor. In 
plan (i) the two-level factors are orthogonal among each other.
In plan (ii) the two and three-level factors are non-orthogonal among each 
other. 

(c)  The two-level factors have C-matrix $(2n)K_2$ in Plan (i) and $nK_2$ in 
Plan (ii). The C-matrices of the three-level factor in Plan (ii) 
 and the four-level factor in Plan  (iii) 
  are given by  $nK_3$ and $nK_4$ respectively. 

(d) The C-matrices of three $n$-level factors A,B and C are as given below.
 $$C_Q = \left ( \left (\begin{array}{cccccc}
2 & -1 & 0 & \cdots & 0 & -1 \end{array} \right )\right ), \; Q
\in G= \{A,B,C\}.$$

 \end{theo}

{\bf Proof : (a) :} Recall the arrays $V$ and $U_i, i=1,2,3 $ described 
in (\ref {2^3array}) and (\ref {U2U3}.
The following arrays represent the plans for the experiments (i), (ii)
 and (iii) respectively. Here the last row represent the block factor.
[Recall Notation \ref {factbl}]. 

$$A_i (4n) = \left [ \begin{array}{cccc}
    U_i  & U_i  & \cdots  &U_i  \\
     V   &  V+1 & \cdots & V + n-1 \\
     \end{array} \right ], i=1,2,3\\$$ 

{\bf (b) :} It can be seen that each of the $n$-level factors satisfies
proportional frequency condition with the two, three and four-level factors. 

{\bf (c) :}
These can be checked by straight computation.

{\bf (d) :} We note the following. Let the $n$-level factors be named as A,B,C in 
that order. Let $G =\{A,B,C\}$. Then, the incidence matrices between them are circulant matrices 
as given below.

\begin{equation}\label{inc-mat}
N^{PQ} = L, P \neq Q, P,Q \in G,\end{equation}

$L$ is as given below.

\begin{equation} \label{L}
 L = ((\begin{array}{cccccc}2 & 1 & 0 & \cdots &0 & 1)) \end{array}.
\end{equation}

Further, one can verify that the pair of factors $P$ and $Q$,
$ P \neq Q, \; P,Q \in G$,  are mutually orthogonal through the
block factor. Now, it
follows from Lemma \ref {orth-relative} that 

$$C_Q = C_{Q;bl}, Q \in G.$$ 

Now, using (\ref {inc-mat}) it is easy to verify that  $C_Q,
Q \in G $ is as given in the statement.   $\Box$ 

{\bf Remark 1:} In the plans constructed in Theorem \ref {constr},
the blocking factor may also be considered as a treatment factor
named D say. Since D is non-orthogonal to each of A,B,C, it's
efficiency would be low. Still it may be useful : particularly for the case 
$n=6$ as using 24 runs we can accomodate
four six-level factors, while on 36 runs one can accomodate at
most three mutually orthogonal six-level factors.

The C-matrix of  D can be obtained 
as follows.
\begin{equation} \label{C-D} C_D = 4I_4 - E_D (H_D)^-(E_D)^T, 
\end{equation} 
where $E_D = \left [ \begin{array}{ccc} M & M & M \end{array} \right ]$
and $H_D = \left [ \begin{array}{ccc}
4I_4 & L    & L \\
L    &  4I_4 & L \\
L    & L     &  4I_4 \\\end{array} \right ].$

Here $M$ is the circulat matrices described below and $L$ is as
in (\ref {L}).
\begin{equation} \label{M,L}
 M = 2((\begin{array}{ccccc}1 & 1 & 0 & \cdots &0)) \end{array} \;
\end{equation}

We now present another construction.

\begin{theo} \label {OMEP-bl} Suppose there exist a connected binary block 
design $d$ with $b$ blocks of size $k$ each and $v$ treatments.
Let $m$ denote the maximum number of constraints of an orthogonal array of 
strength two, $k^2$ runs,  $k$ symbols and index 1.  

(a) Then, $\exists$ an OMEP for a $v^{m-1}$ experiment on $bk$ blocks of size 
$k$ each.

(b) Further, the C-matrix $C_P$ of every factor $P$ is $kC_d$.
\end{theo}

{\bf Proof: } Let us take the $j$ th block of $d$ consisting of the
treatment set $T_j$ say. Let $A_j$ denote the $OA(k^2, m, k,2)$
with $T_j$ as the set of symbols for the first $m-1$ rows and the set
$\{ jk+1, jk+2, \cdots, (j+1)k \}$ for the last row. Then the array

$$A = \left [ \begin{array}{cccc} 
A_1 & A_2 & \cdots & A_b \\ \end{array} \right ].$$

is the mixed array $\rho(bk^2,m; v^{m-1}.bk)$ representing the required 
MEP.

It is easy to verify that condition (\ref {orth-bleq}) is satisfied
by every pair of factor. Thus $A$ represent an OMEP. Part (b)  follows easily.

{\bf Examples :} 1. For a given $v$, there are many connected binary block 
designs with $b$ blocks of size $k$ each. Clearly, it is desirable
that $k$ is a prime or a prime power, so that one can accomodate 
$k$ factors. Again, for the sake of economy, $b$ should not be too large.
Given these conditions, one possibility is to take $b=2$ and $k$
the smallest prime power $\geq [v/2]$. A binary connected design
satisfying the above conditions always exist. Thus, Theorem \ref {OMEP-bl}
implies the existence of the OMEPs for the following experiments.

(i) $4^3$ and $5^3$ on 6 blocks of size 3 each, (ii) $6^4$ and $7^4$ on 8 
blocks of size 4 each, (iii)$8^5$ and $9^5$ on 10 blocks of size 5 each
and so on. 

2. In the examples above the block design $d$ is not equireplicate
and hence the factors would not have equal frequency. To achieve
equal frequency, of course, number of runs has to be bigger, particularly
when $v$ is prime. [ In that case the method of Mukerjee, Dey and Chatterjee 
(2001) might be useful].
We list the following connected and equireplicate block designs
with composite $v$ to be used in Theorem \ref {OMEP-bl} together
with the OMEP obtained.

\vskip5pt

(a) $v=6: b=3, k=4$.
\vspace{.5em} $d : \left[ \begin{tabular}{ccc}
Block 1 & Block 2 & Block 3   \\
\hline \\
1 2 3 4 & 1 2 5 6 & 3 4 5 6 \\
\hline \\
\end{tabular} \right ]$
\vskip5pt

OMEP : $6^4$ experiment on 12 blocks of size 4 each.

(b) $v=8: b= k =4$. 
$ d : \left [ \begin{tabular}{cccc}
Block 1 & Block 2 & Block 3 & Block 4  \\
\hline \\
1 2 3 4 & 5 6 7 8 & 1 2 5 6 & 3 4 7 8 \\
\hline \\
\end{tabular} \right ]$
\vspace{.5em}

OMEP : $8^4$ experiment on 16 blocks of size 4 each.

(c) $v=10 : b=4, k=5$.
$ d : \left [ \begin{tabular}{cccc}
Block 1 & Block 2 & Block 3 & Block 4  \\
\hline \\
1 2 3 4 5 & 6 7 8 9  10 & 1 2 3 6 7 & 4 5 8 9 10 \\
\hline \\
\end{tabular}\right ]$
\vspace{.5em}

OMEP : $10^5$ experiment on 20 blocks of size 5 each.

(d) $v=12 : b=3, k=8$.
$ d : \left [ \begin{tabular}{ccc}
Block 1 & Block 2 & Block 3   \\
\hline \\
1 2 3 4 5 6 7 8 & 1 2 3 4 9 10 11 12 & 5 6 7 8 9 10 11 12\\
\hline \\
\end{tabular}\right ]$
\vspace{.5em}

OMEP : $12^8$ experiment on 24 blocks of size 8 each.

{\bf Remark 3 :} It is interesting to find that the performance
of the OMEPs derived from the block designs (a), (b), (c) and (d) 
above are quite good in the following sense. For each factor in each plan, 
the BLUEs of all the main effect contrsts except two have the same
variance with the hypothetical (real in the case $v=8$) OMEP 
without the blocking factor on the same 
number of runs. The varinces of the BLUEs of the two remaining contrasts
are, of course bigger.

\section{References}
\begin{enumerate}
\item Bagchi,S. (2006). Some series of inter-class orthogonal main
effect plans, Submitted.

\item Bose, M. and Bagchi, S. (2007). Optimal main effect plans on blocks of 
small size,  Jour. Stat. Prob. Let, vol. 77, no. 2, p : 142-147. 

\item Dey, A. and Mukerjee, R.(1999). Fractional factorial
plans, Wiley series in Probability and Statistics, Wiley, New York: .

\item Hedayat,A.S., Sloan, N.J.A. and Stuffken, J. (1999). Orthogonal arrays.
Springer-Verlag, New York.

\item Mukerjee, R., Dey, A. and Chatterjee, K. (2001). Optimal main 
effect plans with non-orthogonal blocking. {\it Biometrika},{\bf 89}, 225-229.

\item Huang, L., Wu, C.F.J. and Yen, C.H. (2002). The idle column method :
Design construction, properties and comparisons, Technometrics, vol.44,
p : 347-368.

\item Ma, C. X., Fang, K.T. and Liski, E (2000). A new approach in 
constructing orthogonal and nearly orthogonal arrays, Metrika,
vol. 50, p : 255-268.

\item Nguyen, (1996). A note on the construction of  near-orthogonal arrays with
mixed levels and economic run size, Technometrics, vol. 38, p : 279-283.

\item Wang, J.C. and Wu.,C.F.J. (1992). Nearly orthogonal arrays with
mixed levels and small runs, Technometrics, vol. 34, p : 409-422.

\item Xu, H (2002). An algorithm for constructing orthogonal and nearly orthogonal 
arrays with mixed levels and small runs, Technometrics, vol. 44, p : 356-368.
\end{enumerate}
****************************************************************88
  \end{document}